\documentclass[11pt]{amsart}

\usepackage{latexsym}
\usepackage{amscd}
\usepackage[all]{xy}
\usepackage[mathscr]{euscript}
\usepackage{dsfont}
\usepackage{hyperref}

\normalsize

\RequirePackage{xspace}

\newcommand{\thought}[1]{}
\renewcommand{\thought}[1]{ \textbf{[#1]}}

\usepackage{enumerate}        
\newenvironment{roenumerate}{\begin{enumerate}[\upshape (i)]}{\end{enumerate}}

\newcommand\nc {\newcommand}
\newcommand\rnc{\renewcommand}

\newcount\blopone
\newcount\xone
\newcount\xtwo
\newcount\ytwo

\newtheorem{theorem}{Theorem}
\newtheorem{prop}[theorem]{Proposition}
\newtheorem{com}[theorem]{Comment}
\newtheorem{apl}[theorem]{Application}
\newtheorem{exercise}[theorem]{Exercise}
\newtheorem{redu}[theorem]{Reduction}
\newtheorem{refinement}[theorem]{Refinement}
\newtheorem{summary}[theorem]{Summary}
\newtheorem{importnota}[theorem]{Important Notation}
\newtheorem{prblm}[theorem]{Problem}
\newtheorem{notation}[theorem]{Notation}
\newtheorem{explanation}[theorem]{Explanation}
\newtheorem{defin}[theorem]{Definition}
\newtheorem{caution}[theorem]{Caution}
\newtheorem{remark}[theorem]{Remark}
\newtheorem{reminder}[theorem]{Reminder}
\newtheorem{illustration}[theorem]{Illustration}
\newtheorem{observation}[theorem]{Observation}
\newtheorem{lemma}[theorem]{Lemma}
\newtheorem{construction}[theorem]{Construction}
\newtheorem{discussion}[theorem]{Discussion}
\newtheorem{corollary}[theorem]{Corollary}
\newtheorem{example}[theorem]{Example}
\newtheorem{conclusion}[theorem]{Conclusion}
\newtheorem{sketch}[theorem]{Sketch}
\newtheorem{triviality}[theorem]{Triviality}
\newtheorem{proto}[theorem]{Prototype Quasifibration}
\newtheorem{cauex}[theorem]{Cautionary Example}
\newtheorem{hypo}[theorem]{Hypothesis}
\newtheorem{subth}{ }[theorem]
\newtheorem{case}{Case}[theorem]
\newtheorem{ssubth}{ }[subth]
\newtheorem{facts}[theorem]{Facts}
\newtheorem{history}[theorem]{Historical Survey}
\newtheorem{proofs}[theorem]{Discussion of the Proofs, Old and New}
\newtheorem{heuristic}[theorem]{Heuristic}

\nc\tri[1]{\begin{triviality}
\label{#1}}
\nc\fac[1]{\begin{facts}
\label{#1}
\begin{em}}
\nc\heu[1]{\begin{heuristic}
\label{#1}
\begin{em}}
\nc\app[1]{\begin{apl}
\label{#1}
\begin{em}}
\nc\skt[1]{\begin{sketch}
\label{#1}
\begin{em}}
\nc\hst[1]{\begin{history}
\label{#1}
\begin{em}}
\nc\pfs[1]{\begin{proofs}
\label{#1}
\begin{em}}
\nc\cas[1]{\begin{case}
\label{#1}
\begin{em}}
\nc\rfn[1]{\begin{refinement}
\label{#1}}
\nc\prt[1]{\begin{proto}
\label{#1}}
\nc\lem[1]{\begin{lemma}
\label{#1}}
\nc\pro[1]{\begin{prop}
\label{#1}}
\nc\thm[1]{\begin{theorem}
\label{#1}}
\nc\dis[1]{\begin{discussion}
\label{#1}
\begin{em}}
\nc\cor[1]{\begin{corollary}
\label{#1}}
\nc\dfn[1]{\begin{defin}
\label{#1}}
\nc\sthm[1]{\begin{subth}
\label{#1}}
\nc\exm[1]{\begin{example}
\label{#1}
\begin{em}}
\nc\obs[1]{\begin{observation}
\label{#1}
\begin{em}}
\nc\plm[1]{\begin{prblm}
\label{#1}
\begin{em}}
\nc\rmk[1]{\begin{remark}
\label{#1}
\begin{em}}
\nc\rmd[1]{\begin{reminder}
\label{#1}
\begin{em}}
\nc\ntn[1]{\begin{notation}
\label{#1}
\begin{em}}
\nc\exe[1]{\begin{exercise}
\label{#1}
\begin{em}}
\nc\xpl[1]{\begin{explanation}
\label{#1}
\begin{em}}
\nc\smr[1]{\begin{summary}
\label{#1}
\begin{em}}
\nc\cau[1]{\begin{caution}
\label{#1}
\begin{em}}
\nc\hyp[1]{\begin{hypo}
\label{#1}}
\nc\imn[1]{\begin{importnota}
\label{#1}
\begin{em}}
\nc\rdn[1]{\begin{redu}
\label{#1}
\begin{em}}
\nc\cax[1]{\begin{cauex}
\label{#1}
\begin{em}}
\nc\cmt[1]{\begin{com}
\label{#1}
\begin{em}}
\nc\con[1]{\begin{construction}
\label{#1}
\begin{em}}
\nc\ill[1]{\begin{illustration}
\label{#1}
\begin{em}}
\nc\ssthm[1]{\begin{ssubth}
\label{#1}
\begin{em}}
\nc\cnc[1]{\begin{conclusion}
\label{#1}
\begin{em}}

\nc\elem{\end{lemma}}
\nc\erdn{\end{em}\end{redu}}
\nc\erfn{\end{refinement}}
\nc\eprt{\end{proto}}
\nc\ethm{\end{theorem}}
\nc\ecor{\end{corollary}}
\nc\edfn{\end{defin}}
\nc\esthm{\end{subth}}
\nc\epro{\end{prop}}
\nc\etri{\end{triviality}}
\nc\eexm{\end{em}
\end{example}}
\nc\eobs{\end{em}
\end{observation}}
\nc\ecmt{\end{em}
\end{com}}
\nc\efac{\end{em}
\end{facts}}
\nc\eheu{\end{em}
\end{heuristic}}
\nc\eapp{\end{em}
\end{apl}}
\nc\ermk{\end{em}
\end{remark}}
\nc\ermd{\end{em}
\end{reminder}}
\nc\eill{\end{em}
\end{illustration}}
\nc\eplm{\end{em}
\end{prblm}}
\nc\ecas{\end{em}
\end{case}}
\nc\eskt{\end{em}
\end{sketch}}
\nc\ecau{\end{em}
\end{caution}}
\nc\ecax{\end{em}
\end{cauex}}
\nc\eimn{\end{em}
\end{importnota}}
\nc\entn{\end{em}
\end{notation}}
\nc\eexe{\end{em}
\end{exercise}}
\nc\expl{\end{em}
\end{explanation}}
\nc\edis{\end{em}
\end{discussion}}
\nc\econ{\end{em}
\end{construction}}
\nc\esmr{\end{em}
\end{summary}}
\nc\ehst{\end{em}
\end{history}}
\nc\epfs{\end{em}
\end{proofs}}
\nc\ehyp{
\end{hypo}}
\nc\ecnc{\end{em}
\end{conclusion}}
\nc\essthm{\end{em}
\end{ssubth}}

\nc\sst{\scriptstyle}
\newcommand{\comment}[1]{}
\newcommand{\ri}{\longrightarrow}
\newcommand{\sr}{\rightarrow}

\newcommand{\zz}{{\mathbb Z}}

\newcommand{\nn}{{\mathbb N}}

\newcommand{\D}{{\mathbf D}}
\newcommand{\qq}{{\mathbb Q}}

\nc\op{^{\hbox{\rm\tiny op}}}
\nc\mth{^{\hbox{\rm\tiny th}}}

\nc\script{\mathscr}
\nc\z{\zeta}
\nc\bc{{\mathbb{BC}}}
\nc\ct{{\script T}}
\nc\cf{{\script F}}
\nc\cg{{\script G}}
\nc\ch{{\script H}}
\nc\ck{{\script K}}
\nc\cl{{\script L}}
\nc\cv{{\script V}}
\nc\ce{{\script E}}
\nc\cs{{\script S}}
\nc\car{{\script R}}
\nc\cd{{\script D}}
\nc\cc{{\script C}}
\nc\ca{{\script A}}
\nc\ci{{\script I}}
\nc\cj{{\script J}}
\nc\co{{\script O}}
\nc\cu{{\script U}}
\nc\cx{{\script X}}
\nc\Cp{{\script P}}
\nc\cq{{\script Q}}
\nc\cy{{\script Y}}
\nc\cz{{\script Z}}
\nc\bd{\begin{description}}
\nc\ed{\end{description}}
\nc\ctob{{\script C}at\big(\ci^{op},\ca\big)}
\nc\clim{{\ds\mathop{\rm lim}_{\ds\longleftarrow}}\,}
\nc\climi{\clim_{\!i}\,}
\nc\climn{\clim^{\!n}\,}
\nc\colim{{\ds\mathop{\rm colim}_{\ds\la}}}
\nc\colimj{{\ds\mathop{\rm colim}_{\ds\la}}{}_{j\,}}
\nc\oa{\overline{\ca}}
\nc\s{\sigma}
\nc\ta{\tau}
\nc\os{\overline\sigma}
\nc\ot{\overline\tau}
\nc\T{\Sigma}
\nc\Tm{\Sigma^{-1}}
\nc\de[1]{{\mathop{\rm deg(#1)}}}
\nc\Ad[1]{\mathop{\rm Ad}(#1)}
\nc\ad[1]{\mathop{\rm ad}(#1)}
\nc\kth{{\it K}--theory}
\nc\loc[1]{{\text{\rm Loc(#1)}}}
\nc\coloc[1]{{\text{\rm Coloc}(#1)}}

\def\der #1 {D\left(#1\right)}
\nc\prf{\begin{proof}}
\nc\eprf{\end{proof}}
\nc\ds{\displaystyle}
\nc\Tor{\text{\rm Tor}}

\nc\cb{{\script B}}
\nc\ab{{\script A}b}

\nc\be{\begin{roenumerate}}
\nc\ee{\end{roenumerate}}

\nc\cat[1]{{\script C}at\Big({\big\{#1\big\}}\op\,\,,\,\,\ab\Big)}
\nc\csab{{\script C}at\big(\cs^{op},\ab\big)}
\nc\ctab{{\script C}at\Big({\{\ct^\alpha\}}^{op},\ab\Big)}
\nc\csex{{\script E}x\big(\cs^{op},\ab\big)}
\nc\ctex{{\script E}x\Big({\{\ct^\alpha\}}^{op},\ab\Big)}
\nc\sub{\qquad\subset\qquad}
\nc\ctr[1]{{\left.\ct\left(-,#1\right)\right|}_{\cs}}
\nc\ctrf[2]{{\left.\ct\left(#1,#2\right)\right|}_{\cs}}
\nc\Ctr[1]{{\left.\ct\left(-,#1\right)\right|}_{\ct^\alpha}}
\nc\Ctrf[2]{{\left.\ct\left(#1,#2\right)\right|}_{\ct^\alpha}}

\nc\la{\longrightarrow}
\nc\nin{\noindent}
\nc\cad[1]{\text{card}(#1)}
\nc\eq{\quad=\quad}
\nc\BA{\begin{array}{c}}
\nc\EA{\end{array}}
\nc\barr{
\[
\begin{array}{cccccccccccccccc}
}
\nc\earr{
\end{array}
\]
}
\nc\as[1]{{\langle S\rangle}^{#1}}
\nc\sh{\text{\it shift}}

\nc\yy[1]{{\left.\ct\left(-,#1\right)\right|}_{\ct^c}}
\nc\vrep[2]{{\left.\ct\left(#1,#2\right)\right|}_{\ct^\alpha}}
\nc\da{\downarrow}
\nc\Hom{{\mathop{\rm Hom}}}
\nc\HHom{{\script H}{\mathop{\rm om}}}
\nc\End{{\mathop{\rm End}}}
\nc\Ext{{\mathop{\rm Ext}}}
\nc\mr\Modtc
\nc\PExt{{\mathop{\rm PExt}}}
\nc\stm{\text{\rm stmod}(kG)}
\nc\stM{\text{\rm StMod}(kG)}
\nc\e{\varepsilon}
\nc\p{\varphi}

\nc\rs{\s^{-1}A}
\nc\br{{\{\s^{-1}A\}}}
\nc\ra\ri
\nc\y[1]{\mathbf{y}#1}
\nc\x[1]{\mathbf{z}#1}
\nc\mmod[1]{#1\text{--\rm mod}}
\nc\Mod[1]{#1\text{--\rm Mod}}
\nc\MMod[1]{\text{\rm Mod--}#1}
\nc\Md {\ensuremath{\mathop{\textup{Mod}}}}
\rnc\mod[1]{\ensuremath{\mathop{#1\textup{--mod}}}\xspace}
\nc\Modtc{\Mod{\ct^c}}
\nc\pgldim[1]{\mathop{\rm pgldim}\,#1}
\nc\tf{{\rm [TR5]}}
\nc\tfs{{\rm [TR5$^*$]}}
\nc\Fun{\text{\rm Funct}(F\op,\ab)}
\nc\sym{\text{\rm Sym}}
\nc\sgn{\text{\rm sgn}}
\nc\Pro{\text{\rm Prod}^{}_\alpha(F\op,\ab)}
\nc\Yt[1]{{\left.\Hom_\ct^{}\left(-,#1\right)\right|}_F^{}}
\nc\dl{\delta}
\nc\Proj[1]{#1\text{--\rm Proj}}
\nc\proj[1]{#1\text{--\rm proj}}
\nc\Flat[1]{#1\text{--\rm Flat}}
\nc\Inj[1]{#1\text{--\rm Inj}}
\nc\Ima{\mathrm{Im}}
\nc\Ker{\mathrm{Ker}}
\nc\ov{\overline}
\nc\wt{\widetilde}
\nc\wh{\widehat}
\nc\ph{\varphi}
\nc\tstr{{\it t}--structure}
\nc\spec[1]{{\text{\rm Spec}(#1)}}

\nc\EProd{\text{\rm EProd}}
\nc\ECoprod{\text{\rm ECoprod}}
\nc\Prod{\text{\rm Prod}}
\nc\ldimp{\text{\rm LDim}^{\prod}}
\nc\ldimc{\text{\rm LDim}^{\coprod}}
\nc\gen[2]{{\langle#1\rangle}^{}_{#2}}
\nc\genu[3]{{\langle#1\rangle}^{[#3]}_{#2}}
\nc\ogen[1]{\ov{\langle#1\rangle}}
\nc\ogenun[2]{\ov{\langle#1\rangle}_{#2}^{}}
\nc\ogenu[3]{\ov{\langle#1\rangle}^{[#3]}_{#2}}
\nc\ogenul[3]{\ov{\langle#1\rangle}^{(-\infty,#3]}_{#2}}
\nc\ogenuf[3]{\ov{\langle#1\rangle}^{[#3,\infty)}_{#2}}
\nc\genuf[3]{{\langle#1\rangle}^{[#3,\infty)}_{#2}}
\nc\genul[3]{{\langle#1\rangle}^{(-\infty,#3]}_{#2}}
\nc\dperf[1]{\D^{\mathrm{perf}}(#1)}
\nc\dcoh{\mathbf{D}^b_{\mathrm{coh}}}

\nc\dmcoh{\mathbf{D}^-_{\mathrm{coh}}}
\nc\dscoh{\mathbf{D}^{}_{\mathrm{coh}}}
\nc\RHHom{{\script{RH}}{\mathrm{om}}}
\nc\Coprod{\mathrm{Coprod}}
\nc\COprod{\mathrm{coprod}}
\nc\add{\mathrm{add}}
\nc\Add{\mathrm{Add}}
\nc\Smr{\mathrm{smd}}
\nc\id{\mathrm{id}}
\nc\LL{\mathbf{L}}
\nc\R{\mathbf{R}}
\nc\wi{\wt{\text{\it\i}}}
\nc\exal{\ce\text{\it x}(\ct^\alpha,\ab)}
\nc\exalz{\ce\text{\it x}_{\aleph_0}^{}(\ct^\alpha,\ab)}
\nc\fc{\mathfrak{C}}
\nc\fl{\mathfrak{L}}
\nc\fs{\mathfrak{S}}
\nc\Prf{\text{\bf Perf}}

\nc\hoco{
\begin{picture}(40,10)
\put(20,0){\makebox(0,0)[b]{\text{\rm Hocolim}}}
\put(5,-2){\vector(1,0){30}}
\end{picture}\,\,}

\nc\holim{
\begin{picture}(40,10)
\put(20,0){\makebox(0,0)[b]{\text{\rm Holim}}}
\put(35,-2){\vector(-1,0){30}}
\end{picture}}

\begin{document}

\author{Amnon Neeman}\thanks{The research was partly supported 
by the Australian Research Council}
\address{Centre for Mathematics and its Applications \\
        Mathematical Sciences Institute\\
        Building 145\\
        The Australian National University\\
        Canberra, ACT 2601\\
        AUSTRALIA}
\email{Amnon.Neeman@anu.edu.au}

\title{Metrics on triangulated categories}

\begin{abstract}
  In a 1973 article Lawvere defined (among many other things)
  \emph{metrics} on categories---the article has been enormously influential
  over the years, spawning
  a huge literature. In recent work, which is surveyed in the current note,
  we pursue a largely-unexplored
  angle: we \emph{complete} categories with respect
  to their Lawvere metrics.

  This turns out to be particularly interesting when the category
  is triangulated and the Lawvere metric is \emph{good;} a metric is
  good if it is translation invariant and the balls of radius $\e>0$
  shrink rapidly enough as $\e$ decreases. The definitions are all
  made precise at the beginning
of the note. And the main theorem is that a certain natural subcategory
$\fs(\cs)$,
of the completion of $\cs$ with respect to a good metric, is triangulated.

There is also a theorem which, under restrictive conditions,
gives a procedure for computing $\fs(\cs)$.
As examples we discuss the special cases 
(1) where $\cs$ is the homotopy category
of finite spectra, and (2) where $\cs=\D^b(\mod R)$, the derived category
of bounded complexes of finitely generated $R$--modules over a noetherian
ring $R$.
\end{abstract}

\subjclass[2010]{Primary 18E30}

\keywords{Triangulated categories, homological functors, metrics}

\maketitle



\rmd{R0.5}
Following a 1973 article of Lawvere~\cite{Lawvere73,Lawvere02}, more
precisely the discussion on
pages~139-140 of~\cite{Lawvere73}\footnote{Lawvere calls
  these \emph{normed categories;} if we wanted to be faithful to
  Lawvere's terminology we would use the word ``norm'' rather than ``metric''.
  But to the author the term ``norm'' suggests that the metric is compatible
  with the action of some ring with a multiplicative absolute value.}, a 
\emph{metric} on a category
is a function that assigns a positive real number (length)
to every morphism, in such a way that for every identity map
$\id:x\la x$ we have $\text{Length}(\id)=0$ and the triangle inequality
is satisfied. The triangle inequality means:
if $x\stackrel f\la y\stackrel g\la z$ are composable morphisms
then
\medskip
\[\text{Length}(gf)\quad\leq\quad\text{Length}(f)+\text{Length}(g)\ .
\]
Lawvere's article does many other things and
has had an enormous influence over the years---when
I last checked on \emph{Google Scholar} it had 732 citations. To the
best of my knowledge the myriad applications have essentially
all gone in directions
totally different from the one we will be pursuing in this note.
There is only a handful of exceptions, we will cite these when they become
relevant.
\ermd

We begin with a string of definitions.

\dfn{D0.3}
Suppose we are given a category $\cc$. Two Lawvere metrics
$\text{\rm Length}_1^{}$ and $\text{\rm Length}_2^{}$ are declared
\emph{equivalent} if,
for any real number $\e>0$, there exists a number
$\delta>0$ such that
\[
\begin{array}{ccc}
  \{\text{\rm Length}_1(a\rightarrow b)<\delta\}&\quad\Longrightarrow\quad&
  \{\text{\rm Length}_2(a\rightarrow b)<\e\}\ , \\
  \{\text{\rm Length}_2^{}(a\rightarrow b)<\delta\}&\quad\Longrightarrow\quad&
  \{\text{\rm Length}_1(a\rightarrow b)<\e\}\ .
\end{array}
\]
\edfn

\dfn{D0.7}
Let $\cc$ be a category with a Lawvere metric.
A \emph{Cauchy sequence} in $\cc$ is a sequence
$E_1\la E_2\la E_3\la\cdots$ of composable morphisms
in which the maps $E_i\la E_j$ eventually become very short.
More precisely: for any $\e>0$ there exists
an $M>0$ such that the morphisms $E_i\la E_j$ satisfy
\[\text{\rm Length}(E_i\sr E_j)\quad<\quad\e
\]
whenever $i,j>M$ and $i\leq j$. 
\edfn

We are accustomed from analysis to the idea of completing a metric
space
with respect
to its metric, and now we want to do the same for Lawvere metrics on
categories.
And the idea is simple enough: the Yoneda embedding takes any category
$\cc$
to a subcategory of a  cocomplete category, with the traditional
definition that a category is
cocomplete if all small colimits exist.
Hence the completion of $\cc$ should just be the
closure of the image under Yoneda of $\cc$. We make this precise in:  

\dfn{D0.9}
Let $\cc$ be a category with a metric. 
Let $Y:\cc\la\Hom[\cc\op,\text{\rm Set}]$ be the Yoneda functor, that
is the functor sending an object $c\in\cc$ to the representable
functor
$Y(c)=\Hom(-,c)$.
\be
\item
Let $\fl'(\cc)$ be the \emph{completion} of $\cc$,
meaning the full subcategory of
$\Hom[\cc\op,\text{\rm Set}]$ whose objects are
the colimits in $\Hom[\cc\op,\text{\rm Set}]$ of Cauchy sequences in $\cc$.
\item
Let $\fc'(\cc)$ be the full subcategory of
$\Hom[\cc\op,\text{\rm Set}]$ whose objects we will call
\emph{compactly supported.} An object 
$F\in\Hom[\cc\op,\text{\rm Set}]$, that is a functor
$F:\cc\op\la\text{\rm Set}$, is declared to be compactly supported if it takes
sufficiently short morphisms to isomorphisms. That is: $F$ belongs to
$\fc'(\cc)$ 
if there
exists an $\e>0$ such that
\[
\{\text{\rm Length}(a\rightarrow b)<\e\}\quad\Longrightarrow\quad\{F(b)\la F(a)\text{ is
  an isomorphism}\}.
\]
\item
Now let $\fs'(\cc)=\fc'(\cc)\cap\fl'(\cc)$.
\ee
\edfn

Next assume the category $\cc$ is \emph{pre-additive.} This means 
that $\Hom(a,b)$ is an abelian group for every pair of objects
$a,b\in\cc$,
and the composition is bilinear.\footnote{In the more recent literature what
  used to be called 
  pre-additive categories goes by the name
  $\zz$--linear categories, or even just
  $\zz$--categories.} In this situation
the Yoneda map factors as a composite
\[\xymatrix@C+30pt{
\cc\ar[r]^-{\wt Y}  &\MMod\cc\ar[r]^-\Phi &\Hom[\cc\op,\text{\rm Set}]
}\]
where $\MMod\cc$ is the
category whose objects are the  \emph{additive functors} of the form
$\cc\op\la\mathrm{Ab}$. And since $\MMod\cc$ is cocomplete we now
have the option of taking the closure of the image of $\wt Y$ instead
of the closure of the image of $Y$. This leads us to

\dfn{D0.10}
Let $\cc$ be a \emph{pre-additive} category with a metric. 
Then

\be
\item
  Let $\fl(\cc)$ be the completion of $\wt Y(\cc)$ in
  $\MMod\cc$; it is the full subcategory of $\MMod\cc$ whose
  objects are the colimits in $\MMod\cc$ of Cauchy sequences in $\cc$.

\item
Let $\fc(\cc)=\Phi^{-1}\fc'(\cc)$.

\item
Finally let $\fs(\cc)=\fc(\cc)\cap\fl(\cc)$.
\ee
\edfn

\rmk{R0.10.5}
We leave it to the reader to compare our description of the completion
with what can be found in Lawvere
\cite[Proposition, bottom of p.~163, and its proof which goes on
to p.~164]{Lawvere73}.

Categories with metrics---that is what Lawvere~\cite{Lawvere73}
calls normed categories---may be viewed
as categories enriched over
a certain closed monoidal category, see Betti and
Galuzzi~\cite{Betti-Galuzzi75} for a detailed exposition. And there
is a notion of completing an enriched category with respect
to a class of colimits, the reader can find it already in
Kelly's book~\cite{Kelly82}, but for a more direct
approach see Kelly and
Schmitt~\cite{Kelly-Schmitt05}. It doesn't seem
automatic that the specialization
of the general theory
to the case at hand agrees with what we've done in this
article. That is: using Kelly's construction
we may complete a normed category
$\cc$,  enriched as in~\cite{Betti-Galuzzi75},
with respect to Cauchy sequences---and what's obtained
doesn't in general seem to agree with our
$\fl'(\cc)$. There are conditions that suffice to guarantee agreement: 
in Kubi{\'s}~\cite{Kubis17} the reader can see that adding the extra
axiom that $\text{\rm Length}(f)\leq\text{\rm Length}(gf)+\text{\rm Length}(g)$
is sufficient. In this article we work 
mostly with  ``good metrics'',
which will
be spelled out in Definition~\ref{D0.15}. Our good metrics
 happen \emph{not} to satisfy
the Kubi\'s axiom. The interested reader can nevertheless check
that restricting to good metrics also suffices to guarantee the
agreement of $\fl'(\cc)$ with the Cauchy completion due to Kelly.

For yet another construction of the category $\fl'(\cc)$ see
Krause~\cite{Krause18}. Krause only looks at one particular metric
but the method generalizes.
In Krause's approach the
category $\fl'(\cc)$ is presented as the Gabriel-Zisman localization
(see~\cite{Gabriel-Zisman67}) of the
category of Cauchy sequences, where one 
formally  inverts the Ind-isomorphisms.

The idea of studying the categories $\fc'(\cc)$ and $\fs'(\cc)$ seems
to have arisen only in \cite{Neeman18A}.
\ermk

\rmk{R0.10.7}
When I wrote \cite{Neeman18A} I was unaware of the earlier work
by Lawvere, Betti, Galuzzi, Kelly, Schmitt and Kubi{\'s};
one of the aims of this survey is to present the
results with the notation as close as possible to the older papers,
but nevertheless
compatible enough with~\cite{Neeman18A} so that the interested reader can
easily read further. Since Lawvere introduces metrics on
arbitrary categories, the right notion of the completion
in his generality
is $\fl'(\cc)$ or Kelly's more sophisticated enriched
completion. In~\cite{Neeman18A} we assume at the outset that
$\cc$ is a triangulated category, hence only mention $\fl(\cc)$.
It is easy to check that, when $\cc$ is pre-additive, the functor
$\Phi:\MMod\cc\la\Hom[\cc\op,\text{\rm Set}]$ restricts
to an equivalence $\fl(\cc)\la\fl'(\cc)$, and hence also to an
equivalence $\fs(\cs)\la\fs'(\cs)$.
\ermk

\rmk{R0.11}
All we have shown so far is that there is no law barring a mathematician
from making a string of ridiculous definitions. To persuade the reader that
this formalism has some value we need to use it to prove a theorem.

In the interest of full disclosure: in the generality of the paragraphs
above I can't prove anything worthwhile. The only obvious observation is
that the constructions are robust under replacing one metric by
an equivalent other. The Cauchy sequences depend only on the equivalence
class of the metric, hence so do the categories $\fl'(\cc)$ and $\fl(\cc)$.
The definitions of the categories $\fc'(\cc)$ and $\fc(\cc)$ make it clear
that these two categories are also unperturbed by replacing a
metric by an equivalent. Hence the same is true for
$\fs'(\cc)=\fc'(\cc)\cap\fl'(\cc)$ and
$\fs(\cc)=\fc(\cc)\cap\fl(\cc)$.

So much for triviality. To get anywhere we need to
narrow our attention considerably.
\ermk

\heu{H0.13}
Let $\cs$ be a \emph{triangulated category.}
We will only consider ``translation invariant''\footnote{The word ``translation'' is sometimes used for the shift functor $\T:\cs\la\cs$; our translation invariance has nothing to do with this translation $\T$.} metrics on $\cs$,
meaning for any homotopy cartesian
square
\[\xymatrix{
a\ar[r]^f\ar[d] & b \ar[d] \\
c\ar[r]^g & d
}\]
we postulate that
\[
\text{Length}(f)\eq\text{Length}(g)\ .
\]
Given any morphism $f:a\la b$ we may form the homotopy cartesian
square  
\[\xymatrix{
a\ar[r]^f\ar[d] & b \ar[d] \\
0\ar[r]^g & x
}\]
and our assumption tells us that 
\[
\text{Length}(f)\eq\text{Length}(g)\ .
\]
Hence it suffices to know the lengths of the morphisms $0\la x$.
Replacing the metric by an equivalent, if necessary, we
may assume our metric takes values in
the set of rational numbers of the form
$\left\{\frac1n\,\,\left|\,\,n\in\nn\right\}\right.$.
To know everything about
the metric it therefore suffices to specify the balls
\[
B_n\eq\left\{x\in\cs\,\,\left|\,\,\text{the morphism }0\la x\text{ has length }\leq\frac1n\right.\right\}\ .
\]
To paraphrase the discussion above: if $f:x\la y$ is a morphism, to compute
its length you complete to a triangle $x\stackrel f\la y\la z$, and then
\[
\text{Length}(f)\eq\inf\left\{\left.\frac1n\,\,\,\right|\,\, z\in B_n\right\}\ .
\]
Furthermore we will restrict our attention to 
non-archimedian metrics, that is metrics that satisfy the strong triangle
inequality. This means:
if $x\stackrel f\la y\stackrel g\la z$ are composable morphisms,
then
$\text{Length}(gf)\leq\max\big(\text{Length}(f),\text{Length}(g)\big)$.
By the translation-invariance it suffices to consider the
case $x=0$; that is it suffices to show that the composable morphisms
$0\stackrel f\la y\stackrel g\la z$
satisfy $\text{Length}(gf)\leq\max\big(\text{Length}(f),\text{Length}(g)\big)$.
Completing $g$ to a triangle $y\stackrel g\la z\la w$ this comes down
to $\{y,w\in B_n\}\Longrightarrow\{z\in B_n\}$.
\eheu

The discussion above motivates

\dfn{D0.15}
Let $\cs$ be a triangulated category. A \emph{good metric} on $\cs$ is a
sequence of full subcategories $\{B_n,\,n\in\mathbb{N}\}$, 
containing $0$ and with $B_1=\cs$, and furthermore satisfying
\be
\item
$B_n*B_n=B_n$, which means that if there
exists a triangle
$b\la x\la b'$ with $b,b'\in B_n$, then $x\in B_n$.
\item
$\Tm B_{n+1}\cup B_{n+1}\cup \T B_{n+1}\subset B_{n}$.
\ee
\edfn

\rmk{R0.17}
In Heuristic~\ref{H0.13} we explained where part (i) of
Definition~\ref{D0.15} comes from, it guarantees that the translation-invariant
metric given by the balls $B_n$ is non-archimedian. The hypothesis (ii) of
Definition~\ref{D0.15} has not yet been motivated. We clearly must have
$B_{n+1}\subset B_n$, the ball of radius $\frac 1{n+1}$ must be contained in
the ball of radius $\frac1n$. But it turns out to be convenient to assume the
balls decrease rapidly enough for the stronger hypothesis (ii) to hold; it
guarantees that the automorphism $\T$ is a ``homeomorphism'' with respect to
the metric---in other words the metric $\{\T B_n,\,n\in\nn\}$ is equivalent to
the metric $\{B_n,\,n\in\nn\}$.

Note that we are \emph{not} assuming that the metric is compatible with
any other automorphism of $\cs$.
\ermk

\exm{E0.19}
Suppose $\cs$ is a triangulated category, $\ca$ is an abelian category
and $H:\cs\la\ca$ is a homological functor. Put $B_1=\cs$. If for
$n>1$ we set $B_n$ as given in the formulas
below, we obtain three (inequivalent) good metrics on $\cs$.
\be
\item
$B_n=\{s\in\cs\mid H^i(s)=0\text{ for all $i$ in the range }i>-n\}$.
\item
$B_n=\{s\in\cs\mid H^i(s)=0\text{ for all $i$ in the range }i<n\}$.
\item
$B_n=\{s\in\cs\mid H^i(s)=0\text{ for all $i$ in the range }-n< i< n\}$
\ee
Note that if $\{B_n,\,n\in\nn\}$ define a good metric on $\cs$ then
$\{B\op_n,\,n\in\nn\}$ define a good metric on $\cs\op$, which we will
call the dual metric. Now a homological functor $H:\cs\la\ca$ has a dual
$H\op:\cs\op\la\ca\op$, and the reader can check
that (i), applied to $H\op:\cs\op\la\ca\op$, gives a good metric
equal to the dual of that obtained from
(ii) applied to $H:\cs\la\ca$.

The metric of (iii)
is self-dual.
\eexm

One more definition before the first theorem.

\dfn{D0.21}
Let $\cs$ be a triangulated category with a good metric. With the
category $\fs(\cs)$ as in Definition~\ref{D0.10}(iii), we define the
distinguished  triangles
in  
$\fs(\cs)$ to be the colimits in 
$\fs(\cs)\subset\text{\rm Mod--}\cs$ of Cauchy sequences
of distinguished triangles in $\cs$.
\edfn

\xpl{X0.23}
What this means is the following. If we are given a Cauchy
sequence of distinguished triangles in $\cs$, we can always form the
colimit in the cocomplete category $\MMod\cs$,
and by the definition of $\fl(\cs)$ this colimit must lie in $\fl(\cs)$.
In general there is no guarantee that the colimit will lie in the
subcategory $\fs(\cs)\subset\fl(\cs)$. What Definition~\ref{D0.21} does
is declare that those colimits which happen to
lie in $\fs(\cs)$ are distinguished
triangles in $\fs(\cs)$.
\expl

And now we come to

\thm{T0.25}
With the 
distinguished triangles as in Definition~\ref{D0.21}, the category
$\mathfrak{S}(\cs)$ is triangulated.
\ethm

\nin
The proof of a slightly stronger theorem [the hypotheses on
  the metric are slightly less restrictive]
may be found in \cite[Theorem~2.11]{Neeman18A}.

\rmk{R0.27}
Up to Theorem~\ref{T0.25} all we saw was a string of increasingly bizarre
definitions. We've said it before: in this free world of ours
there is no
law prohibiting
a mathematician from making up a long sequence of
absurd-looking definitions.

Then, out of all the seemingly pointless formalism, we magically pulled
out Theorem~\ref{T0.25}. Perhaps it takes an expert to appreciate how
surprising the result is. Triangulated categories have been around since
the early 1960s---meaning for about 55 years. And the conventional wisdom
has always been that they don't reproduce. Until very recently there were
no interesting recipes that began with a triangulated category $\cs$, and
out of it cooked up another triangulated category $\ct$---in this context
we view
(full) triangulated subcategories and Verdier quotients as
dull, trivial constructions. In more detail: the
first whiff of such a recipe came in 2005 in Keller~\cite{Keller05}.
Keller proved that, given a
triangulated category $\cs$ and an automorphism $\s:\cs\la\cs$, then the
category $\ct=\cs/\s$ sometimes [rarely] has a triangulated structure so that
the quotient map is triangulated---but the conditions are very stringent.
And the only other known
recipe
was found in 2011 by Balmer~\cite{Balmer11A}: given a separable monoid $R$
in a tensor triangulated category
$\cs$, the category $\ct$ whose objects are the $R$--modules in $\cs$,
and whose morphisms are those morphisms in $\cs$ which respect the
$R$--module structure,
is triangulated. The distinguished triangles  in $\ct$ are
precisely those sequences $T\la T'\la T''\la\T T$ 
whose image in $\cs$ is a distinguished triangle.
OK: these two relatively recent exceptions aside, the
accepted wisdom has long been that you need some enhancement to produce
triangulated categories in a nontrivial way.

Theorem~\ref{T0.25} gives a third recipe, and the natural question is whether
the end product is of any value. Given an input triangulated category $\cs$,
together with its good metric, is the output
category $\fs(\cs)$ a worthwhile object
of study?
And to answer this we need examples.
In Example~\ref{E0.19} we saw three ways to
produce good metrics, out of any homological functor $H:\cs\la\ca$.
For each of these the question
arises: what is the triangulated category
$\fs(\cs)$? Is it of any interest?

In general I don't know how to compute $\fs(\cs)$. The only procedure
I know so far assumes that $\cs$ has an embedding into a
larger triangulated
category $\ct$, and this embedding satisfies a strong condition.
In the presence of such an embedding $\fs(\cs)$ may be
computed as a triangulated subcategory of $\ct$. Below we will spell out
carefully the
exact statements.
\ermk

In order to make the above precise we need some more definitions---my
apologies to the reader, we will get to the point in Theorem~\ref{T0.33}.

\dfn{D0.29}
Let $\cs$ be a triangulated category with a
good metric.
Suppose we are given a
fully faithful triangulated functor $F:\cs\la\ct$; we consider
also the functor $\cy:\ct\la\MMod\cs$, which takes an object $A\in\ct$
to the functor $\Hom\big(F(-),A\big)$. The functor $F$ is called a
\emph{good extension with respect to the metric} if $\ct$ has
countable coproducts,
and for every Cauchy sequence $E_*$ in $\cs$ the natural map
$\colim\,\wt Y(E_*)\la \cy\big(\hoco F(E_*)\big)$ is an isomorphism.
\edfn

\xpl{X0.30}
The functors $F$, $\cy$ and $\wt Y$ are related by a canonical
natural isomorphism $\wt Y\cong\cy\circ F$.
And we remind the reader: given a sequence $T_1\la T_2\la T_3\la\cdots$
of composable morphisms in $\ct$, the homotopy colimit is
defined to be the third edge of the triangle
\[\xymatrix@C+20pt{
\ds\coprod_{i=1}^\infty T_i \ar[rr]^-{\id\,-\,\,\text{\rm shift}} & &
\ds\coprod_{i=1}^\infty T_i \ar[r]^-\ph &
\hoco T_*
}\]
where $\text{\rm (shift)}:\coprod_{i=1}^\infty T_i\la\coprod_{i=1}^\infty T_i$ is
the unique map rendering commutative, for any integer $n\geq1$, the square below
\[\xymatrix@C+20pt{
T_n\ar[rr]\ar[d]_{\text{\rm (inc)}_n}& &T_{n+1}\ar[d]^{\text{\rm (inc)}_{n+1}}\\ 
\ds\coprod_{i=1}^\infty T_i \ar@{.>}[rr]^-{\text{\rm shift}} & &
\ds\coprod_{i=1}^\infty T_i 
}\]
with $(\text{\rm inc})$ being the canonical inclusion into the coproduct.
The object $\hoco T_*$ is only defined up to (non-canonical)
isomorphism in $\ct$, but
the isomorphism can be assumed to
respect the map $\ph$. Hence any such isomorphism, between two candidates
for $\hoco T_*$, will respect all the composites
$T_n\stackrel{\text{\rm inc}}\la\coprod_{i=1}^\infty T_i\stackrel\ph\la\hoco T_*$;
we write these as $\ph_n:T_n\la \hoco T_*$. The vanishing of the
composite $\ph\circ(\id-\text{\rm shift})$, in the displayed
maps of the triangle above, guarantees that, for each
integer $n\geq1$, the composite $T_n\la T_{n+1}\stackrel{\ph_{n+1}}\la\hoco T_*$
must be equal to $\ph_n:T_n\la\hoco T_*$.

If we are given a sequence $E_1\la E_2\la E_3\la\cdots$ in the category
$\cs$, then the functor $F$ takes it to a sequence
$F(E_1)\la F(E_2)\la F(E_3)\la\cdots$ in the category $\ct$.
The paragraph above gives, for each $n$, a 
map $\ph_n:F(E_n)\la \hoco F(E_*)$. Applying to this the functor $\cy$
we deduce the second morphism in the composable pair below
\[\xymatrix@C+20pt{
\wt Y(E_n) \ar[r]^-\sim & \cy F(E_n)\ar[r]^-{\cy(\ph_n)} &
\cy\big(\hoco F(E_i)\big)
}\]
where the first morphism
comes from the canonical isomorphism $\wt Y\cong \cy F$.
And these maps assemble to a single morphism
$\colim\,\wt Y(E_*)\la \cy\big(\hoco F(E_*)\big)$, unique up
to (non-canonical) isomorphism. Hence 
postulating that this map is an isomorphism, as in Definition~\ref{D0.29},
makes sense independent of choices.

Note that this is a strong restriction. If $\cs$ has countable coproducts
we might be tempted to let $F$ be the identity $\id:\cs\la\cs$. But then
it becomes a strong hypothesis to assume that the Cauchy sequences
all satisfy the condition that $\colim\,\wt Y(E_*)\la \wt Y\big(\hoco E_*\big)$
is an isomorphism.
\expl

\dfn{D0.31}
Suppose $\cs$ is a triangulated category with
a good metric, and let 
$F:\cs\la\ct$ be a good extension. We define
\be
\item
The full subcategory  $\wh\fl(\cs)\subset\ct$ has for objects
all the homotopy colimits of Cauchy sequences in $\cs$.
\item
The full subcategory $\wh\fs(\cs)\subset\ct$ is
given by the formula
$\wh\fs(\cs)=\wh\fl(\cs)\cap\cy^{-1}\big(\fc(\cs)\big)$, with
$\fc(\cs)\subset\MMod\cs$ as in Definition~\ref{D0.10}(ii).
\ee
\edfn

\thm{T0.33}
The category $\wh\fs(\cs)$ is a
triangulated subcategory of $\ct$, and the functor
$\cy:\ct\la\MMod\cs$ restricts to a
triangulated equivalence $\cy:\wh\fs(\cs)\la\fs(\cs)$.
\ethm

\nin
The proof of a slightly stronger theorem [once again, the hypotheses on
  the metric are slightly less restrictive]
may be found in \cite[Theorem~3.15]{Neeman18A}.

\rmk{R0.999}
We have a fully faithful functor $\wt Y:\cs\la \MMod\cs$ and,
in the presence of a good extension, another
fully faithful functor 
$F:\cs\la\ct$. If we confuse $\cs$ with its essential images
we can view it as a subcategory in each of $\MMod\cs$ and $\ct$.
And then we have subcategories $\wh\fs(\cs)\subset\ct$ and
$\fs(\cs)\subset\MMod\cs$, and
it's natural to wonder what one can say about the subcategories
$\cs\cap\wh\fs(\cs)\subset\ct$ and $\cs\cap\fs(\cs)\subset\MMod\cs$.
To avoid getting too confused, between the incarnation of $\cs$
as a subcategory of $\ct$ and as a subcategory of $\MMod\cs$,
for most of this remark our notation will be careful; we will not confound
$\cs$ with either of its images.

The functor $F:\cs\la\ct$ is a fully faithful, triangulated functor,
while the subcategory $\wh\fs(\cs)\subset\ct$ is triangulated.
Hence  $F^{-1}\big[\wh\fs(\cs)\big]$ is a triangulated subcategory of $\cs$,
and the functor $F$ restricts 
to a fully faithful, triangulated functor
$F^{-1}\big[\wh\fs(\cs)\big]\la\wh\fs(\cs)$.

Now $\wh\fs(\cs)$ lies in $\wh\fl(\cs)\subset\ct$, and
the functor $\cy:\ct\la\MMod\cs$ obviously takes $\wh\fl(\cs)\subset\ct$
to $\fl(\cs)\subset\MMod\cs$; hence $\cy$
restricts to a functor
$\cy|^{}_{\wh\fl(\cs)}:\wh\fl(\cs)\la\fl(\cs)$.
And it turns out to be easy to show that the functor
$\cy|^{}_{\wh\fl(\cs)}$ is essentially surjective, full and conservative.
This means:
every object in $\fl(\cs)$ is isomorphic to an object in the
image of $\cy|^{}_{\wh\fl(\cs)}$, any morphism in $\fl(\cs)$ between objects
in the image of $\cy|^{}_{\wh\fl(\cs)}$ is in the image of $\cy|^{}_{\wh\fl(\cs)}$,
and a morphism in $\wh\fl(\cs)$ is an isomorphism if and only if
$\cy|^{}_{\wh\fl(\cs)}$ takes it to an isomorphism. And the relevance of this
for us is that the commutative square
\[\xymatrix{
\wh\fs(\cs) \ar@{^{(}->}[r]\ar[d]_{\cy|^{}_{\wh\fs(\cs)}} & 
  \wh\fl(\cs)\ar[d]^{\cy|^{}_{\wh\fl(\cs)}}\\
\fs(\cs) \ar@{^{(}->}[r] & \fl(\cs)
}\]
is a strict pullback square. The point is that, from their definitions,
the categories $\wh\fs(\cs)$ and $\fs(\cs)$ are replete subcategories
in, respectively, $\ct$ and $\MMod\cs$; this means they contain all isomorphs
of any of their objects. Theorem~\ref{T0.33} tells us that
the vertical map on the left is an equivalence---hence any object
$x\in\fs(\cs)$ is isomorphic to $\cy(z)$  with 
$z$ an object of $\wh\fs(\cs)$. But if we have an object $t\in\wh\fl(\cs)$
with $\cy(t)=x\cong\cy(z)$, then the isomorphism must lift to
$\wh\fl(\cs)$, and hence $t\cong z$ must belong to the replete
subcategory $\wh\fs(\cs)\subset\wh\fl(\cs)$.

Now recall the Yoneda embedding $\wt Y:\cs\la\fl(\cs)\subset\MMod\cs$.
We have the triangle of functors
\[\xymatrix@C+50pt@R-20pt{
 & \wh\fl(\cs)\ar[dd]^{\cy|^{}_{\wh\fl(\cs)}}\\
\cs\ar[ur]^-F\ar[dr]_-{\wt Y} & \\
 & \fl(\cs)
}\]
which commutes up to natural isomorphism.
It immediately follows that
\[
\wt Y^{-1}\big[\fs(\cs)\big]\eq
F^{-1} \cy^{-1}_{\wh\fl(\cs)}\big[\fs(\cs)\big]\eq
F^{-1}\big[\wh\fs(\cs)\big]
\]
The first equality is because the inverse images
under the isomorphic functors
$\wt Y\simeq F\circ\big[\cy|^{}_{\wh\fl(\cs)}\big]$, of the
replete subcategory $\fs(\cs)$, must be equal.
And the second equality comes from the paragraph above, which informs
us that $\cy^{-1}_{\wh\fl(\cs)}\big[\fs(\cs)\big]=\wh\fs(\cs)$.

Rewriting the second paragraph of the current Remark,
by appealing to the equality
$\wt Y^{-1}\big[\fs(\cs)\big]=F^{-1}\big[\wh\fs(\cs)\big]$,
we deduce first that
$\wt Y^{-1}\big[\fs(\cs)\big]$ is a triangulated subcategory of
$\cs$, and then that the functor $\wt Y:\cs\la\MMod\cs$
restricts to a fully faithful, triangulated functor
$\wt Y^{-1}\big[\fs(\cs)\big]\la\fs(\cs)$.

All of the discussion above assumed we were in the presence of
a good extension $F:\cs\la\ct$. But the assertion
of the last paragraph turns out to
be robust. Even though the category $\fl(\cs)$ is rarely triangulated
it contains both $\cs$ and $\fs(\cs)$ as subcategories---in
the case of $\fs(\cs)$ this is by definition, while for
$\cs$ we commit the notational crime of
confusing $\cs$ with its its essential image under $\wt Y:\cs\la\fl(\cs)$.
Each of $\cs$ and $\fs(\cs)$ has its
own triangulated structure. And it is always true that
$\cs\cap\fs(\cs)$ has a (unique) triangulated structure so that
each of
the two embeddings, into $\cs$ and into $\fs(\cs)$, is triangulated.
\ermk

The import of Theorem~\ref{T0.33} and Remark~\ref{R0.999}
is that any good extension of
$\cs$ contains both $\cs$ and $\fs(\cs)$ as triangulated subcategories,
and the embedding of $\fs(\cs)$ into $\ct$ is explicit enough to facilitate
computations, both
of $\fs(\cs)$ and of $\cs\cap\fs(\cs)$.
The author will be the first to admit that better computational
tools would be wonderful---this is all we have right now.

Notwithstanding the current limitations on what we know, Theorem~\ref{T0.33}
does produce interesting examples. We give a few.

\exm{E0.35}
Let $\cs$ be the homotopy category of finite spectra.
Let us remind the reader: the objects in this category
may be taken to be pairs $(X,n)$, where $X$ is
a pointed, finite CW-complex and $n\in\zz$ is an integer,
positive or negative---the way to
think of this is that the object $(X,n)$ is the $n\mth$ suspension
of $X$.  And $\T:\cs\la\cs$ is the functor
taking a pair $(X,n)$ to the pair $(\T X,n)$, where $\T X$ is the ordinary
suspension of the pointed CW-complex $X$.
The morphisms
are precisely what one would expect, given that we want to
force the functor $\T:\cs\la\cs$ to be invertible---for
any two objects $(X,m)$ and $(Y,n)$ in $\cs$,
the abelian group $\Hom_\cs^{}\big[(X,m),(Y,n)\big]$ is defined to
be the
colimit as $k\la\infty$
of the (eventual) abelian groups $\Hom_{\text{CW-complexes}}^{}(\T^{m+k}X,\T^{n+k}Y)$.
This means: if $k$ is large enough, so that both $m+k$ and $n+k$ are
$\geq 2$, then the $\Hom$-set above is the abelian group of homotopy equivalence
classes of pointed continuous maps $\T^{m+k}X\la\T^{n+k}Y$.

Let
$H:\cs\la\text{\rm Ab}$ be the homological functor which takes a 
spectrum $(X,n)$ to its zeroth stable homotopy group;
in the notation above this means
\[H(X,n)\eq\colim\,\Hom_{\text{CW-complexes}}^{}(\mathbb{S}^k,\T^{n+k}X)\ ,\]
where $\mathbb{S}^k$ is the $k$--dimensional sphere.
In the standard notation of homotopy theorists
$H(X,n)=\pi_0^{}(\T^nX)$ and $H^i(X)=\pi_{-i}^{}(\T^nX)$, where $\pi_{-i}^{}$
is the $(-i)\mth$ stable homotopy group.
Now let the good metric be as
in Example~\ref{E0.19}(i).

Let $F:\cs\la\ct$ be the embedding
of the homotopy category of finite spectra into the
homotopy category of all spectra---the homotopy category of
all spectra is not quite so easy to describe simply, hence
let us leave this out. For us what's important is that the
functor $F$ can be shown to be a good extension. And the computation of
$\wh\fs(\cs)$, which by Theorem~\ref{T0.33} is canonically triangle
equivalent to
$\fs(\cs)$, can be carried out. It shows that $\wh\fs(\cs)\subset\ct$ is
given by the formula
\[
\wh\fs(\cs)\eq\left\{x\in\ct\left|
\begin{array}{c}
H^i(x)=0\text{ for all but finitely many }i\in\zz,\text{ and}\\
H^i(x)\text{ is a finitely generated $\zz$--module for all }i\in\zz
\end{array}
\right.\right\}\ .
\]

The assertions in the paragraphs above follow from the far
more general \cite[Example~4.2]{Neeman18A}.

With $\cs$ still as above, it's known that
$\cs\cap\wh\fs(\cs)\cong\cs\cap\fs(\cs)=\{0\}$
and that $\cs$ and $\fs(\cs)$ are not triangle equivalent.
Let us recall.

Let $\mathbb{S}\in\cs$ be the zero-sphere;  in the notation of the
first
paragraph of the current Example this means 
$\mathbb{S}=(\mathbb{S}^n,-n)$ where $n>0$ is an integer and
$\mathbb{S}^n$ is the $n$--dimensional sphere. And let $K(\zz,0)$
be the 
Eilenberg-MacLane spectrum defined by 
\[
\pi_{-i}^{}\big[K(\zz,0)\big]\eq H^i\big[K(\zz,0)\big]\eq
\left\{
\begin{array}{lll} \zz &\quad&\text{if }i=0\\
0 & & \text{otherwise}
\end{array}
\right.
\]
Next we adopt the terminology of Bondal and 
Van den Bergh~\cite[2.1]{BondalvandenBergh04}: an object $G$ in a
triangulated category $\car$ is a \emph{classical generator} if the 
smallest thick subcategory containing $G$ is all of $\car$.
It's not difficult to see that $\mathbb{S}$ is a classical generator
for $\cs$ while $K(\zz,0)$ is a classical generator for $\wh\fs(\cs)$.
But it's also known that
\[
\begin{array}{cll} 
 \Hom(\mathbb{S},\T^n\mathbb{S})=0 &\quad&\text{if }0<n\\
 \Hom\big(K(\zz,0),\T^n K(\zz,0)\big)=0& & \text{if }n<0\\
\end{array}
\]
and that, for all $n\in\zz$, the modules above
are finitely generated $\zz$--modules.
The vanishing assertions are ancient, and the finite
generation was proved in Serre's 1951 PhD thesis~\cite{Serre51}.
Somewhat more recent is the computation that
\[
\begin{array}{cll} 
\Hom\big(K(\zz,0),\T^n\mathbb{S}\big)=0& & \text{for all }n\neq1
\end{array}
\]
while $\Hom\big(K(\zz,0),\T\mathbb{S}\big)$ is a $\qq$--vector space:
the reader can find this in Lin's 1976 article~\cite[Theorem 3.6]{Lin76}.
The interested reader can look
at  Margolis 1974 article~\cite{Margolis74} for 
an (independent) approach to results similar to Lin's,
and at 
Ravenel~\cite[Section~4]{Ravenel84} for later developments
and extensions.
Anyway: because $\mathbb{S}$ and $K(\zz,0)$ 
classically generate the respective subcategories, we immediately
deduce
\[
\begin{array}{cclcll} 
\text{For any pair }G,H\in\cs&\quad& \Hom(G,\T^n H)=0 &\quad&\text{ if }&0\ll n\\
\text{For any pair }G,H\in\wh\fs(\cs)&\quad& \Hom\big(G,\T^n H\big)=0&\quad & \text{ if }&n\ll0 \\
\text{For any pair }G\in\wh\fs(\cs),\,H\in\cs&\quad& \Hom\big(G,\T^n H\big)=0&\quad & \text{ if }&|n|\gg0
\end{array}
\]
while for all $n$ we have 
\[
\begin{array}{cclcll} 
\text{For any pair }G,H\in\cs&\quad& \Hom(G,\T^n H)&\quad&\text{ is a f.g.~$\zz$--module}\\
\text{For any pair }G,H\in\wh\fs(\cs)&\quad& \Hom\big(G,\T^n H\big)&\quad &\text{ is a f.g.~$\zz$--module}  \\
\text{For any pair }G\in\wh\fs(\cs),\,H\in\cs&\quad& \Hom\big(G,\T^n H\big)&\quad &\text{ is a $\qq$--vector space}
\end{array}
\]

What's more it's known that $\Hom(G,\T^nG)\neq0$ for
infinitely many $n$, when $G$ is either $\mathbb{S}$ or
$K(\zz,0)$; these 
estimates on the non-vanishing of Hom-sets can also be found 
in
Serre~\cite{Serre51}.
It immediately follows that the categories
$\cs$ and $\wh\fs(\cs)$
cannot be triangle equivalent, and that
$\cs\cap\wh\fs(\cs)\cong\cs\cap\fs(\cs)=\{0\}$.
\eexm

\exm{E0.898}
A different example comes about as follows. Let $R$ be a noetherian ring,
let $\D^b(\mod R)$ be the derived category whose objects are the bounded
complexes of finitely generated $R$--modules, and let
$H:\D^b(\mod R)\la(\mod R)$ be the homological functor taking an object
of $\D^b(\mod R)$ to its zeroth cohomology module.
We take on $\big[\D^b(\mod R)\big]\op$
the good metric given by applying
Example~\ref{E0.19}(ii) to $H\op$.

Let $\D(\Mod R)$ be the unbounded derived category of all complexes of
$R$--modules. Then the natural inclusion
$F:\big[\D^b(\mod R)\big]\op\la\big[\D(\Mod R)\big]\op$ is a good extension
with respect to the metric. And this can be used to compute
$\wh\fs\Big(\big[\D^b(\mod R)\big]\op\Big)$ as a subcategory
of $\big[\D(\Mod R)\big]\op$. It turns out to be
$\big[\ch^0\big(\Prf(R)\big)\big]\op$, where $\ch^0\big(\Prf(R)\big)$ is the
derived category whose objects are bounded
complexes of finitely generated, projective $R$--modules.
In this particular case there is an inclusion---and
Remark~\ref{R0.999} tells us that, with $\cs=\D^b(\mod R)$,
the
subcategories $\cs\op$ and $\fs(\cs\op)$
of $\fl(\cs\op)$ satisfy
$\fs(\cs\op)\subset\cs\op$. Furthermore this inclusion respects the
triangulated structure.

The assertions in the paragraph above follow from the far
more general \cite[Proposition~5.6]{Neeman18A}.

Thus out of the category $\D^b(\mod R)$ we have cooked
up its triangulated subcategory $\ch^0\big(\Prf(R)\big)$, and hence we also
know the quotient
\[
\D_{\text{\rm sing}}^{}(R)\eq\frac{\D^b(\mod R)}{\ch^0\big(\Prf(R)\big)} \ ,
\]
where $\D_{\text{\rm sing}}^{}(R)$ is what's known in the literature
as the singularity category of $R$.

OK: the paragraphs above showed that, given the category $\D^b(\mod R)$
and its metric, then out of the data we can construct the
triangulated subcategory
$\ch^0\big(\Prf(R)\big)$ and the quotient $\D_{\text{\rm sing}}^{}(R)$. The
reader might naturally ask if there is a way to construct the metric
without appealing to the homological functor $H:\D^b(\mod R)\la(\mod R)$.
The answer turns out to be Yes up to equivalence. The equivalence
class of the metric can be obtained
without using anything other than the triangulated structure on
$\D^b(\mod R)$---we will say a little more about this in Remark~\ref{R???}. 

It is also possible to construct examples where the inclusion goes
the other way, that is $\cs\subset\fs(\cs)$.
In fact: starting with the triangulated category
$\cs=\ch^0\big(\Prf(R)\big)$ and the homological functor $H:\ch^0\big(\Prf(R)\big)\la(\mod R)$,
the functor taking a cochain complex to its zeroth cohomology module,
we can endow $\cs$ with the good metric of Example~\ref{E0.19}(i).
And then it may be computed
that $\fs(\cs)=\D^b(\mod R)$, and as subcategories of
$\fl(\cs)$ we have an inclusion $\cs\subset\fs(\cs)$
which agrees with the standard triangulated inclusion
$\ch^0\big(\Prf(R)\big)\subset\D^b(\mod R)$. Once again: the 
equivalence class of the metric has an intrinsic description, it
depends only on the triangulated structure of $\cs=\ch^0\big(\Prf(R)\big)$.
Thus $\ch^0\big(\Prf(R)\big)$ also contains
enough data to determine $\D_{\text{\rm sing}}^{}(R)$. See Remark~\ref{R???} 
for some elaboration, and~\cite{Neeman18A} for infinitely more detail
and much greater generality.
\eexm

\rmk{R???} 
In Example~\ref{E0.898} we asserted that, up to equivalence, certain metrics
have intrinsic descriptions. It's time to explain this.

Let $\cs$ be a triangulated category and let $G\in\cs$ be a classical
generator---recall: classical generators were defined in 
Bondal and 
Van den Bergh~\cite[2.1]{BondalvandenBergh04}, and we have already met them
in Example~\ref{E0.35}. The functor $\Hom(G,-)$ is a
homological functor, and we may apply the formulas
of Example~\ref{E0.19} to obtain metrics. It is an easy exercise
to show that, if $G,G'\in\cs$ are two classical generators,
then the metrics
obtained by applying Example~\ref{E0.19}~(i), (ii) or (iii) to
the functor $\Hom(G,-)$ are (respectively) equivalent to the metrics obtained
by applying Example~\ref{E0.19}~(i), (ii) or (iii) to
the functor $\Hom(G',-)$. Thus on any triangulated category
with a classical generator there are three intrinsic 
equivalence classes of metrics. For the
sake of definiteness let us call them 
$\cl\!\text{\it ength}_{(i)}^{\mathrm{Int}}$
$\cl\!\text{\it ength}_{(ii)}^{\mathrm{Int}}$ and 
$\cl\!\text{\it ength}_{(iii)}^{\mathrm{Int}}$. 
To spell it out explicitly: to obtain
a metric in the equivalence class $\cl\!\text{\it ength}_{(i)}^{\mathrm{Int}}$
you choose a classical generator $G$ for the category $\cs$, and then
the formula of Example~\ref{E0.19}(i), applied to the homological
functor $\Hom(G,-)$, delivers your metric.

The category $\ch^0\big(\Prf(R)\big)$ is known to have a classical 
generator, more precisely the object
$R\in\ch^0\big(\Prf(R)\big)$ is a
classical generator. 
To obtain a metric in the equivalence class 
$\cl\!\text{\it ength}_{(i)}^{\mathrm{Int}}$
we begin by choosing $R$ for our classical generator,
and then we apply the
formula of Example~\ref{E0.19}(i) to the homological functor $\Hom(R,-)$.
But there is a canonical 
isomorphism $H(-)=\Hom(R,-)$, between of the functor $H$ of final paragraph of
Example~\ref{E0.898} and the functor $\Hom(R,-)$.
Hence the metric of
the final paragraph of
Example~\ref{E0.898} is in the equivalence
class $\cl\!\text{\it ength}_{(i)}^{\mathrm{Int}}$, it is intrinsic up
to equivalence.

We also asserted that the metric on $\big[\D^b(\mod R)\big]\op$, 
which was studied
in the opening paragraphs of Example~\ref{E0.898}, is intrinsic up
to equivalence. Here the argument is subtler, partly because
(as far as I know) the category $\D^b(\mod R)$ does not have to have a 
classical generator. The main idea on how to define the metric 
intrinsically is presented in the opening paragraphs of the introduction
to~\cite{Neeman18A}.

There are cases in which the category $\D^b(\mod R)$ is known to have a 
classical generator, and moreover this classical generator is known
to satisfy some additional nice properties: it is ``strong''---whatever
that means---and one has further 
technical knowledge about it. When this 
happens to be the case the intrinsic description of the metric 
simplifies. Concretely: if the ring $R$ is commutative, noetherian,
has finite Krull-dimension and every closed subset of $\spec R$ admits
a regular alteration, then $\D^b(\mod R)$ has a
classical generator, and the metric on
$\big[\D^b(\mod R)\big]\op$ given in the first paragraph of 
Example~\ref{E0.898} is in the equivalence
class $\cl\!\text{\it ength}_{(ii)}^{\mathrm{Int}}$.
\ermk

\exm{E0.7934}
In Remark~\ref{R0.10.5} we briefly mentioned the article~\cite{Krause18} by 
Krause.  Let us discuss his work a little more fully.

As in the article \cite{Neeman18A}, Krause~\cite{Krause18} comes up with a 
procedure that
produces out of $\cs$ the category we call 
$\fl(\cs)$---Krause's recipe is different 
from Definition~\ref{D0.10}(i), but out of a different oven 
comes exactly the same dish. For the reader interested in looking up
Krause~\cite{Krause18} for more detail: the category
we call $\fl(\cs)$ is canonically equivalent to
what goes by the name  $\wh\cs$ in 
\cite{Krause18}. And just as in 
\cite{Neeman18A}  Krause
looks at the special case
where $\cs=\ch^0\big(\Prf(R)\big)$ as in Example~\ref{E0.898}.
With $H:\ch^0\big(\Prf(R)\big)\la (\mod R)$ the homological
functor of the last paragraphs of Example~\ref{E0.898}, Krause studies
the metric of Example~\ref{E0.19}(iii)---this is where
his treatment radically differs from
Example~\ref{E0.898}, where our metric was the one of Example~\ref{E0.19}(i).
Because we now have two metrics let us denote the completion
of Example~\ref{E0.898} by $\fl_1(\cs)$ and Krause's 
completion by $\fl_2(\cs)$. It isn't difficult to show that, inside
the category $\MMod\cs$, there is an inclusion $\fl_1(\cs)\subset\fl_2(\cs)$.

Of course the category $\fs_1(\cs)\cong\D^b(\mod R)$, 
being a full subcategory of $\fl_1(\cs)$, is
also a full subcategory of $\fl_2(\cs)$. And there is an intrinsic 
description of the subcategory $\fs_1(\cs)\subset\fl_2(\cs)$,
the reader can find it in~\cite{Krause18}. Krause denotes it $\wh\cs^b$,
I suppose in our notation it should be $\fl_2(\cs)^b$.

Next we come to the triangles. With respect to the metric
of Example~\ref{E0.19}(i) the triangles in $\fs_1(\cs)$ are simply
the colimits of Cauchy sequences of triangles in $\cs$---subject
of course to the restriction that the colimit lies in
$\fs_1(\cs)\subset\fl_1(\cs)$. 
Theorem~\ref{T0.25} tells us that, with the triangles defined as above, the
category $\fs_1(\cs)$ is triangulated.

Now let's work out some 
consequences. 
Let $F'\la F$ be any morphism in 
$\fs_1(\cs)$; it can be expressed as the colimit of a Cauchy sequence
of morphisms $s'_*\la s_*$ in $\cs$, where the metric on $\cs$ is as
in Example~\ref{E0.19}(i). We may (non-canonically) complete this to 
a sequence of triangles  $s'_*\la s_*\la s''_*\la \T s'_*$ in $\cs$, 
and the reader can easily check that this sequence is Cauchy
in the metric of Example~\ref{E0.19}(i) and the colimit lies
in $\fs_1(\cs)$. By Definition~\ref{D0.21} the colimit
is a distinguished triangle
in $\fs_1(\cs)$ extending the morphism $F'\la F$. But from the axioms
of triangulated categories the extension of the morphism $F'\la F$ to
a distinguished triangle in $\fs_1(\cs)$ is unique up to non-canonical
isomorphism. It follows that the sequences 
of triangles in $\cs$ of the form $s'_*\la s_*\la s''_*\la \T s'_*$,
extending the given Cauchy sequence of morphisms $s'_*\la s_*$,
must all be non-canonically Ind-isomorphic.

This can be proved, but the proof I know relies heavily on the fact that 
the colimit lies in
the subcategory $\fs_1(\cs)\subset\fl_1(\cs)$, which
 was chosen carefully in terms
of the metric. Whereas Krause's definition of $\fl_2(\cs)^b\subset\fl_2(\cs)$
makes no mention of the metric.

Now let's compare the Cauchy sequences with respect to the
two metrics under consideration. Since the metric of 
Example~\ref{E0.19}(i) is finer than the metric of 
Example~\ref{E0.19}(iii) there are more Cauchy sequences with respect 
to Krause's metric---this is what leads to the (proper) inclusion
$\fl_1(\cs)\subset\fl_2(\cs)$. As it turns out any Cauchy
sequence of morphisms $s'_*\la s_*$, 
with respect to the metric of Example~\ref{E0.19}(iii) and whose
colimit happens to lie in $\fl_1(\cs)\subset\fl_2(\cs)$, is Ind-isomorphic
to a Cauchy sequence $t'_*\la t_*$ with respect to the metric
of Example~\ref{E0.19}(i). So we might be tempted to guess that
the Cauchy sequences of triangles with respect to the metric
of Example~\ref{E0.19}(iii), with colimits in $\fs_1(\cs)=\fl_2(\cs)^b$,
will also be Ind-isomorphic to Cauchy sequences of triangles
with respect to the metric of Example~\ref{E0.19}(i).
But if we try to produce
such an Ind-isomorphism we
run  into the problem that the mapping cone isn't functorial---the 
simple-minded  
approach breaks down. As the definition
of $\fl_2(\cs)^b\subset\fl_2(\cs)$ doesn't involve the metric
I see no sophisticated alternative to the simple-minded 
method---for all I know
there might be Cauchy sequences of triangles, with respect to
the metric of Example~\ref{E0.19}(iii) and with
colimit in $\fs_1(\cs)$, which aren't Ind-isomorphic to Cauchy sequences
of triangles with respect to
the metric of Example~\ref{E0.19}(i).

Krause's solution to the problem is to fix  an enhancement,
and only admit those Cauchy sequences
that lift to the chosen enhancement. For the situation 
at hand a minimal enhancement
suffices---it's enough to assume we are working with Keller's towers,
see Keller~\cite{Keller91} for the original exposition, or his
appendix to Krause~\cite{Krause18} for a condensed version.

Of course it is possible to apply the machinery surveyed here
to Krause's metric---we obtain a triangulated category $\fs_2(\cs)$ whose
triangulated structure is enhancement-free. 
Using a good extension with respect to the metric and 
Theorem~\ref{T0.33} it can be computed
that $\fs_2(\cs)$
is a proper subcategory of $\fs_1(\cs)\cong\D^b(\mod R)$---the objects
are those complexes in $\D^b(\mod R)$ which have bounded
injective resolutions. For more detail
the reader is referred to \cite[Example~4.9]{Neeman18A}.
\eexm

We have said it before but repeat for emphasis: Theorem~\ref{T0.33} is
at present the only computational tool we have. It would be great to have some
more ways to compute $\fs(\cs)$.

\rmk{R0.3077}
Krause's metric can be defined on other triangulated categories, for example
on the homotopy category of finite spectra of Example~\ref{E0.35}.
Thus the category of finite spectra also has two completions
$\fl_1(\cs)\subset\fl_2(\cs)$, the one
corresponding to the metric of Example~\ref{E0.35}
being contained in the one with respect to Krause's 
metric. In Barthel's appendix to 
Krause's paper~\cite[Appendix~A]{Krause18} there is a proof that,
just as in the case of the $\cs$ studied in Example~\ref{E0.7934},
for this $\cs$ too we have the equality $\fs_1(\cs)=\fl_2(\cs)^b$. 
\ermk

\rmk{R0.37}
In this survey we've tried to convince the reader that good metrics on
triangulated categories can be useful. We've only touched on what's
possible---the reader interested in more theorems in
this vein is referred to the
longer and more extensive survey~\cite{Neeman17B}. To give 
one instance of
a result
in~\cite{Neeman17B} which
is immediately relevant to our discussion above:
Examples~\ref{E0.35}
and \ref{E0.898} may look quite different, but both can be obtained as
special cases of a single, much more general
example.
\ermk

\medskip

\nin{\bf Acknowledgements.} \ \ 
The author would like to thank Henning Krause for asking the question that
led to the research we have surveyed in this
note, which is expounded more fully and in greater generality 
in~\cite{Neeman18A}. The author also wants to thank
Steve Lack for pointing him to Lawvere's old
paper~\cite{Lawvere73,Lawvere02}, Richard Garner for explaining the 
relevance of Betti and Galuzzi~\cite{Betti-Galuzzi75}, 
Kelly~\cite{Kelly82}, Kelly and Schmitt~\cite{Kelly-Schmitt05} and 
Kubi\'s~\cite{Kubis17}, as well as Paul Balmer, Ivo Dell'Ambrogio,
David Gepner, Bernhard Keller, 
Sasha Kuznetsov,
Tom Leinster, Nori Minami, 
Doug Ravenel and Greg Stevenson for corrections and helpful comments
on earlier drafts of this manuscript and help with references to other 
relevant literature.

\def\cprime{$'$}


\begin{thebibliography}{BVdB03}

\bibitem[Bal11]{Balmer11A}
Paul Balmer.
\newblock Separability and triangulated categories.
\newblock {\em Adv. Math.}, 226(5):4352--4372, 2011.

\bibitem[BG75]{Betti-Galuzzi75}
Renato Betti and Massimo Galuzzi.
\newblock Categorie normate.
\newblock {\em Boll. Un. Mat. Ital. (4)}, 11(1):66--75, 1975.

\bibitem[BVdB03]{BondalvandenBergh04}
Alexei~I. Bondal and Michel Van~den Bergh.
\newblock Generators and representability of functors in commutative and
  noncommutative geometry.
\newblock {\em Mosc. Math. J.}, 3(1):1--36, 258, 2003.

\bibitem[GZ67]{Gabriel-Zisman67}
Peter Gabriel and Michel Zisman.
\newblock {\em Calculus of fractions and homotopy theory}.
\newblock Ergebnisse der Mathematik und ihrer Grenzgebiete, Band 35.
  Springer-Verlag New York, Inc., New York, 1967.

\bibitem[Kel82]{Kelly82}
G.~Maxwell Kelly.
\newblock {\em Basic concepts of enriched category theory}, volume~64 of {\em
  London Mathematical Society Lecture Note Series}.
\newblock Cambridge University Press, Cambridge-New York, 1982.

\bibitem[Kel91]{Keller91}
Bernhard Keller.
\newblock Derived categories and universal problems.
\newblock {\em Comm. in Algebra}, 19:699--747, 1991.

\bibitem[Kel05]{Keller05}
Bernhard Keller.
\newblock On triangulated orbit categories.
\newblock {\em Doc. Math.}, 10:551--581, 2005.

\bibitem[Kra]{Krause18}
Henning Krause.
\newblock Completing perfect complexes.
\newblock https://arxiv.org/abs/1805.10751.

\bibitem[KS05]{Kelly-Schmitt05}
G.~Maxwell Kelly and Vincent Schmitt.
\newblock Notes on enriched categories with colimits of some class.
\newblock {\em Theory Appl. Categ.}, 14:no. 17, 399--423, 2005.

\bibitem[Kub]{Kubis17}
Wies{\l}aw Kubi\'s.
\newblock Categories with norms.
\newblock https://arxiv.org/abs/1705.10189.

\bibitem[Law73]{Lawvere73}
F.~William Lawvere.
\newblock Metric spaces, generalized logic, and closed categories.
\newblock {\em Rend. Sem. Mat. Fis. Milano}, 43:135--166 (1974), 1973.

\bibitem[Law02]{Lawvere02}
F.~William Lawvere.
\newblock Metric spaces, generalized logic, and closed categories [{R}end.
  {S}em. {M}at. {F}is. {M}ilano {\bf 43} (1973), 135--166 (1974); {MR}0352214
  (50 \#4701)].
\newblock {\em Repr. Theory Appl. Categ.}, (1):1--37, 2002.
\newblock With an author commentary: Enriched categories in the logic of
  geometry and analysis.

\bibitem[Lin76]{Lin76}
Tsau~Young Lin.
\newblock Duality and {E}ilenberg-{M}ac{L}ane spectra.
\newblock {\em Proc. Amer. Math. Soc.}, 56:291--299, 1976.

\bibitem[Mar74]{Margolis74}
Harvey~R. Margolis.
\newblock Eilenberg-{M}ac{L}ane spectra.
\newblock {\em Proc. Amer. Math. Soc.}, 43:409--415, 1974.

\bibitem[Nee21]{Neeman17B}
Amnon Neeman.
\newblock Approximable triangulated categories.
\newblock In {\em Representations of Algebras, Geometry and Physics}, volume
  769 of {\em Contemp. Math.}, pages 111--155. Amer. Math. Soc., Providence,
  RI, 2021.

\bibitem[Neeb]{Neeman18A}
Amnon Neeman.
\newblock The categories $\ct^c$ and $\ct^b_c$ determine each other.
\newblock https://arxiv.org/abs/1806.06471.

\bibitem[Rav84]{Ravenel84}
Douglas~C. Ravenel.
\newblock Localization with respect to certain periodic homology theories.
\newblock {\em Amer. J. of Math.}, 106:351--414, 1984.

\bibitem[Ser51]{Serre51}
Jean-Pierre Serre.
\newblock Homologie singuli\`ere des espaces fibr\'{e}s. {A}pplications.
\newblock {\em Ann. of Math. (2)}, 54:425--505, 1951.

\end{thebibliography}
\end{document}